\documentclass{amsart}

\usepackage{amsmath}
\usepackage{amsfonts}
\usepackage{amsthm}
 
\newcommand{\fsmpsp}{\left(\Omega, \mathcal{F}, \{\mathcal{F}_t\}, P \right)}
\newcommand{\e}{\mathrm{E}}
\newcommand{\mcal}[1]{\mathcal{#1}}

\newcommand{\ltwo}{{\mathbf{L}^2}}

\newcommand{\norm}[1]{\left\lVert #1 \right\rVert}

\newcommand{\abs}[1]{\left\lvert #1\right\rvert}

\newcommand{\rr}{\mathbb{R}}

\renewcommand{\mcal}[1]{\mathcal{#1}}

\newcommand{\eq}{\begin{equation}}
\newcommand{\en}{\end{equation}}

\newcommand{\re}[1]{\mbox{(\ref{#1})}}
\newcommand{\rem}[1]{\mbox{\em (\ref{#1})}}

\newcommand{\wien}{\mathbb{W}}
\newcommand{\sort}{\mcal{S}}
\newcommand{\qvar}[1]{\langle #1 \rangle}

\newcommand{\qmu}{Q\cdot\mu}
\newcommand{\pmu}{P\cdot\mu}
\newcommand{\qmun}[1]{Q_{#1}\cdot\mu_{#1}}
\newcommand{\qjmun}[2]{Q_{#1}\cdot\mu_{#2}}

\newcommand{\qjmu}[1]{Q_{#1}\cdot\mu}

\newcommand{\combi}[2]{\begin{pmatrix}#1 \\ #2 \end{pmatrix}}

\newtheorem{Theorem}{Theorem}
\newtheorem{theorem}[Theorem]{Theorem}
\newtheorem{lemma}[Theorem]{Lemma}

\newtheorem{corollary}[Theorem]{Corollary}
\newtheorem{construction}[Theorem]{Construction}
\newtheorem{proposition}[Theorem]{Proposition}
\newtheorem{example}[Theorem]{Example}
\newtheorem{exercise}[Theorem]{Exercise}
\newtheorem{defn}[Theorem]{Definition}
\newtheorem{claim}[Theorem]{Claim}
\newtheorem{question}[Theorem]{Question}
\newtheorem{conjecture}[Theorem]{Conjecture}
\newtheorem{condition}[Theorem]{Condition}
\newtheorem{remark}[Theorem]{Remark}
\newtheorem{problemma}[Theorem]{Problemma}
\newtheorem{jpfigure}[Theorem]{Figure}
\newtheorem{jptable}[Theorem]{Table}

\theoremstyle{definition}

\def\proof{\noindent{\bf Proof.\ \ }}

\def\endpf{\hfill $\Box$ \vskip .25in}

\newfont{\msbm}{msbm10 at 12pt}
\newfont{\eusb}{eusb10}
\newfont{\eusm}{eusm10}
\newfont{\eurb}{eurb10}
\newfont{\eurm}{eurm10}
\newfont{\eufb}{eufb10}
\newfont{\eufm}{eufm10}

\newcommand {\ints} {\mbox{\msbm\symbol{'132}}}

\newcommand {\reals} {\mathbb{R}}
\newcommand {\PR} {\mathbb{P}}

\newcommand {\hf} { \mbox{$ {1 \over 2 }$} }

\newcommand{\te}{\rightarrow}
\newcommand{\ed}{\mbox{$ \ \stackrel{d}{=}$ }}
\newcommand{\convd}{\mbox{$ \ \stackrel{\!d}{\rightarrow}$ }}
\newcommand{\eps}{\varepsilon}

\newcommand{\comment}[1]{}
\newcommand{\noshowcomment}{\renewcommand{\comment}[1]{}}

\newcommand{\lb}[1]{\label{#1}}

\newenvironment{cnj}[1]{\begin{conjecture}\label{#1}}{\end{conjecture}}

\newenvironment{crl}[1]{\begin{corollary}\protect\label{#1}}{\end{corollary}}

\noshowcomment
\noshowcomment

\newcommand{\la}{\lambda}
\newcommand{\de}{\delta}
\newcommand{\De}[1] { \Delta_{#1} }

\newcommand{\FF}{{\mcal F}}

\begin{document}

\title[Brownian motions interacting through ranks]{One-dimensional Brownian particle systems with rank dependent drifts}

\author{Soumik Pal}
\address{506 Malott Hall\\ Cornell University\\ Ithaca NY 14853}
\email{soumik@math.cornell.edu}

\author{Jim Pitman}
\address{367 Evans Hall \# 3860\\
University of California\\
Berkeley CA 94720}
\email{pitman@stat.berkeley.edu}


\thanks{Soumik Pal's research is partially supported by N.S.F. grant DMS-0306194 to the probability group at Cornell.}
\thanks{Jim Pitman's research is supported in part by N.S.F. Grants DMS-9970901, DMS-0405779, and DMS-0071448.}

\date{\today}
\maketitle
\begin{abstract}We study interacting systems of linear Brownian motions whose drift vector at every time point is determined by the relative ranks of the coordinate processes at that time. Our main objective has been to study the long range behavior of the spacings between the Brownian motions arranged in increasing order. For finitely many Brownian motions interacting in this manner, we characterize drifts for which the family of laws of the vector of spacings is tight, and show its convergence to a unique stationary joint distribution given by independent exponential distributions with varying means. We also study one particular countably infinite system, where only the minimum Brownian particle gets a constant upward drift, and prove that independent and identically distributed exponential spacings remain stationary under the dynamics of such a process. Some related conjectures in this direction have also been discussed. 
\end{abstract}

\section{Introduction}
In this paper we consider systems of interacting one-dimensional Brownian motions $X = (X_i(t), i\in I,\; t \ge 0)$, where $i$ ranges over an index $I$, which is either the finite set $\{1, \ldots, N\}$, or the countable set of positive integers $\mathbb{N}$. 

For $I=\{1,2,\ldots,N\}$, if we define ordered coordinates of any vector $x\in\rr^N$ by
\[
x_{(1)} \le x_{(2)} \le \ldots \le x_{(N)},
\]
the locations $X_i(t)$ of the Brownian particles evolve according to the system of stochastic differential equations
\eq
\lb{sde1}
d X_i (t) = \sum_{j\in I} \delta_j 1\left\{\;X_i(t)= X_{(j)}(t)\;\right\} dt + dB_i(t)   ~~~~~~~~(i \in I)
\en
for some sequence of {\em drifts} $\delta_1, \delta_2, \ldots \in \reals$.
Here the $B_i$'s are asumed to be independent $(\FF_t)$-Brownian motions for some suitable underlying filtration $(\FF_t)$. 
Less formally, the Brownian particles evolve independently except that the $i$th ranked particle is given drift $\delta_i$.

For these finite systems, with arbitrary initial distribution of $(X_i(0), i \in I)$, and arbitrary drifts $\delta_i$, the existence and uniqueness in law of such an $N$ particle system is guaranteed by the standard theory of 
SDE's. See Lemma \ref{lm2} for more details. We would also like to define SDE \eqref{sde1} when $I$ is the entire countably infinite collection $\mathbb{N}$. This is trickier since the ordered process might not remain well-defined anymore, and the existence of the solution of SDE \eqref{sde1} depends both on the initial distribution and the sequence of drifts. We consider one such system in Section \ref{sec.lag} with a drift sequence $(\delta, 0, 0, \ldots)$ ($\delta > 0$) for which the finite-dimensional arguments can be suitably extended. For this system to exist in the weak sense, starting from an initial vector $X_1(0) < X_2(0) < \ldots$, we show that it suffices to assume 
\[
\liminf_{n\rightarrow \infty} \frac{(X_n(0)-X_1(0))^2}{n} >0.
\]
See Lemma \ref{infchange} for the details of the proof.

Consider now the {\em Brownian spacings system} derived from this ordered Brownian particle system with rank dependent drifts, i.e.,
$$
\De{k}(t):= X_{(k+1)}(t) - X_{(k)}(t)  \mbox{ for } k, k + 1  \in I. 
$$

The ordered particle system derived from independent Brownian motions with no drift (meaning $\delta_i \equiv 0$) was studied by Harris \cite{harris65}, Arratia \cite{arratia83}
and Sznitman \cite{sznitman86},\cite[p. 187]{sznitman91}.
By Donsker's theorem, this system can be interpreted as a scaling limit of ordered particle systems derived from independent symmetric nearest neighbour random walks on $\ints$. Harris \cite{harris65} considers the spacings of an infinite ordered Brownian particle system defined by : 
$$
\Delta_i^*(t) := B_{(i+1)}(t) - B_{(i)}(t) ~~~~~~~~(i \in \ints)
$$
where $\{B_i\}$ is a family of independent Brownian motions with no drifts
and initial states $B_i(0) = B_{(i)}(0)$ such that $B_0(0) = 0$ and the
$B_i(0)$ for $i \in \ints \backslash\{0\}$ are points of a Poisson process of rate $\la$
on $\reals$. That is, $(B_{(i)}(t), t \ge 0 )_{i \in \ints}$ is the almost surely
unique collection of processes with continous paths such that $B_{(i)} (t) \le B_{(i+1)}(t)$ for all
$i \in \ints, t \ge 0$ and the union of the graphs of these processes is identical to the
union of the graphs of the Brownian paths $(B_{i}(t), t \ge 0 )_{i \in \ints}$.
We call $(B_{(i)}(t), i \in \ints, t \ge 0 )$ the {\em Harris system of ordered Brownian motions with rate $\la$}
and their differences $(\Delta_{i}^*(t), i \in \ints, t \ge 0 )$ the 
{\em Harris system of Brownian spacings with rate $\la$}.

As observed by Arratia \cite[\S 4]{arratia85},
\comment{(according to Ferrari-Fontes this goes back to Kesten: check Spitzer 1970 \cite{xxx}.)}
the corresponding stationary system of spacings between particles of the exclusion process associated with a nearest-neighbor random walk 
on $\ints$ can be interpreted as a finite or infinite series of queues, also known as the {\em zero range process with constant rate}.
See \cite{harrison73,harrison00,harvan97,harrison78,harwil87s,ocyor01} for background on systems of Brownian queues. Such connections between systems of queues and one-dimensional interacting
particle systems have been exploited by a number of authors,
in particular Kipnis \cite{kipnis86}, 
Srinivasan \cite{srinivasan93}, Ferrari-Fontes \cite{ferfo94}, \cite{ferf96}, and Sepp\"{a}l\"{a}inen \cite{sepp97}.
Ferrari \cite{fer96} surveys old and new results on the limiting behavior of a tagged particle in various interacting particle systems. Also see articles by Baryshnikov \cite{barysh01} and O'Connell-Yor \cite{ocyor01} for some recent studies of Brownian queues in tandem, connected to the directed percolation and the directed polymer models, and the GUE random matrix ensemble. 

More recently rank dependent SDEs have been considered by several authors as possible models for financial or economic data. Fernholz in \cite{fern02} introduces the so-called {\em Atlas model} which we study in this paper. It is a model of finitely many Brownian particles, where at every time point the minimum Brownian motion gets a constant positive drift, while the rest of them get no drift. The general rank dependent interacting Brownian models, whose drifts and volatilities depend on time-varying ranks have been considered by Banner, Fernholz, and Karatzas in \cite{atlasmodel}, with whom our work in this paper bears close resemblance. For SDE \eqref{sde1}, they work under a specific condition on the drift sequence required for stability of the solution process. We prove in  Theorem \ref{theoremN} (condition \eqref{conal}) that this condition is indeed necessary and sufficient. Although their method is mostly based on a beautiful analysis of the local times of intersections of different Brownian motions, they also note the connection with the Harrison-Williams theory of reflected Brownian motions which we use extensively in this paper. See Lemma \ref{lm1} for a complete statement. 

In that same article \cite{atlasmodel} the authors establish marginal convergence of spacings $\Delta_k$ to exponential distributions. They leave the question of joint convergence open, which we settle in this paper in Theorem \ref{theoremN} by proving that the vector of spacings converge jointly to independent exponentials with different means. They also study ergodic properties of such processes including a demonstration of the exchangeability of the indices of the Brownian particles under rank-dependent drifts and volatilities. In their later papers, Fernholz and Karatzas also consider generalized versions of \eqref{sde1} where the drifts and the volatilities depend on both the index and the rank of a Brownian particle. As expected, in most such cases, explicit descriptions of their properties become very difficult to obtain. However, many interesting results can still be recovered. A good source for what has been done so far can be found in the recent survey article by Fernholz and Karatzas \cite{karsurv}. Also see the article by Chatterjee and Pal \cite{chatpal07} for a follow-up in this direction where the authors consider an increasing number of Brownian particles in a rank-dependent motion and establish connection with the Poisson-Dirichlet family of point processes.

McKean and Shepp, in \cite{sheppmckean}, also consider Brownian motions interacting via their ranks. They start with two Brownian motions, and their objective is to find the optimal drift under constraints (as a control) such that the probability that both Brownian motions never hit zero is maximized. The solutions is, as they establish, the Atlas model for the two particle system.

An interesting related model, studied by Rost and Vares in \cite{rostvares}, replaces the ordered particles in the Harris model by linear Brownian motions repelled by their nearest neighbors through a potential. The authors study stationary measures for the spacings of such processes, and show that rescaled combinations of spacings converge to an Ornstein-Uhlenbeck process. 

Our purpose here is to draw attention to the general class of Brownian particle systems with rank dependent drifts, as considered in \cite{atlasmodel}. Many natural questions about these systems remain open. We are particularly interested in an infinite version of the Atlas model, with a drift sequence $(\delta,0,0, ...)$ with \emph{Atlas drift} $\delta >0$. One result we obtain for this system is the following:

\begin{theorem}\label{thpois}
For each $\de >0$,
a sequence of independent Exponential$(2 \de)$ variables provides an 
equilibrium distribution of spacings for the infinite Atlas model with the Atlas drift $\delta$.
\end{theorem}

Theorem \ref{thpois}, which essentially follows from Theorem \ref{derivatzero}, suggests a number of interesting conjectures and 
open problems. In particular, we can immediately formulate:

\begin{cnj}{cnjpois1}
For each $\de >0$, Theorem \ref{thpois} describes the
unique equilibrium distribution of spacings for the infinite Atlas particle system 
with Atlas drift $\delta$.
\end{cnj}

Let
$(\De{1}(t), \De{2}(t), \ldots )_{t \ge 0}$ 
denote the equilibrium state 
of spacings of the infinite Atlas Brownian particle 
system described by Theorem \ref{thpois}.
This process has some subtle features.
For each $k = 1,2, \ldots$ and each $t >0$, 
\eq
\lb{def}
(\De{k}(t), \De{k+1}(t), \ldots ) \ed (\De{1}(t), \De{2}(t), \ldots )
\en
and the common distribution of these sequences is that of independent
exponential $(2 \de)$ variables.
But while both sides of \re{def} define stationary 
sequence-valued processes as $t$ varies, these
processes do not have the same law for all $k$.
In particular, the finite-dimensional
distributions of non-negative stationary process 
$(\De{k}(t), t \ge 0)$
depend on $k$.

Harris \cite[eqn. (7.1)]{harris65} gave an explicit formula for the law of
$B_{(0)}(t)$, the location at time $t$ of particle initially at
$0$ in the Harris system of ordered Brownian motions,
from which he deduced for $2 \lambda = 1$ that 
\eq\label{har0}
{B_{(0)}(t) \over  t ^{1/4} } \convd 
\left(2 \over \pi \right)^{1/4} {1 \over \sqrt{2 \lambda}} \, B(1) \mbox{ as } t \te \infty 
\en
where $B(1)$ is standard Gaussian.
As remarked by Arratia \cite[p. 71]{arratia83}, the conclusion for general
$\lambda >0$ folows from the case $2 \lambda = 1$ by Brownian scaling.
See also De Masi, Ferrari \cite{demasiferrari}, Rost, Vares \cite{rostvares}, and Arratia \cite{arratia83} where variants (or generalizations) of \re{har0} is proved for a tagged particle in the exclusion process on $\ints$ associated with 
a simple symmetric random walk.
Harris conjectured that the process $B_{(0)}$ is not Markov, and 
left open the problem of describing the long-run
behaviour of paths of $B_{(0)}$. These questions were answered by D\"urr, Goldstein, and Lebowitz in \cite[Thm 7.1]{dgl85} where they prove that the \emph{tagged} process $B_{(0)}(t)$, suitably rescaled, converges to fractional Brownian motion with Hurst parameter $1/4$. In fact, they show a general theorem where convergence to fractional Brownian motion holds for systems with processes which have stationary increments and perform \emph{elastic collisions} as in the Harris model. Also see \cite{dgl87} where the same authors generalize their results in the case where an external potential is present.

In this connection, we have the following conjecture:
\begin{cnj}{cnjhar}
For each $\de >0$, in the infinite Atlas model with Atlas drift $\delta >0$, and initially ordered values $X_{(k)}(0),k = 1,2, \ldots$ which are the points of a Poisson process with rate $2 \de $ on $(0,\infty)$,
$$
{ X_{(k)}(t) - X_{(k)} (0) \over  t ^{1/4} } \convd 
\left(2 \over \pi\right)^{1/4} {c_k \over \sqrt{2 \de}} \, B(1)  \mbox{ as } t \te \infty
$$
for some sequence of constants $c_k >0$ with $c_k \te 1$ as $k \te \infty$.
\end{cnj}

The paper is organized as follows. In the next section, we analyze the finite Brownian particle system with rank dependent drifts. The main result is Theorem \ref{theoremN}, which describes the convergence in total variation of the laws of the spacings to that of independent exponential distributions. The precise condition needed on the drift sequence for such stability has also been proved. In Section \ref{sec.lag} we look at countably infinite Brownian particles with the dynamics of the Atlas model. The main result, Theorem \ref{thpois}, follows readily from Theorem \ref{derivatzero}.

\section{The finite Brownian particle system}
\label{sec.fin}

We first present some results regarding the asymptotic behavior of spacings for the $N$ particle system with
arbitrary rank-dependent drifts $\delta_i, 1 \le i \le N$.
We start by recording the following two lemmas, which clarify the issues of
existence and uniqueness of the $N$ particle system with arbitrary drifts $\delta_i$,
and characterize the associated ordered particle system.
See the book by Revuz and Yor \cite{ry99} for a background and the definitions of concepts from the calculus of
continuous semimartingales.

\begin{lemma}\label{lm1} 
Let $X_i, 1 \le i \le N$ be a solution of the SDE \re{sde1} defined on some filtered probability space $\fsmpsp$, for some arbitrary initial condition and arbitrary drifts $\{\delta_j\}$.
Then for each $1 \le j \le N$ the $j$th ordered process $X_{(j)}$ is a continuous semimartingale 
with decomposition
\begin{equation}\label{tanak}
dX_{(j)}(t) = d \beta_{j}(t) + 
{ 1 \over \sqrt{2}} ( d L_{(j-1,j)}(t) -  d L_{(j,j+1)}(t)  )
\end{equation}
where the $\beta_{j}$'s for $1 \le j \le N$ are independent $(\FF_t)$-Brownian motions with
unit variance coefficient and drift coefficients $\delta_j$, and
where $L_{(0,1)} = L_{(N,N+1)} =  0,$
and for $1 \le j \le N-1$ 
\eq
\lb{occd}
L_{(j,j+1)}(t) = \lim_{ \eps \downarrow 0 } {1 \over 2 \eps }  \int_0^t 1 ( (X_{(j+1)} (s)- X_{(j)}(s))/\sqrt{2} \le \eps ) ds
~~~(t \ge 0)
\en
which is half the continuous increasing local time process at $0$ of the semimartingale 
$(X_{(j+1)} - X_{(j)})/\sqrt{2}$. Moreover,
the ordered system is a Brownian motion in the domain
\eq
\lb{domn}
\left\{ (x_{(i)}, 1 \le i \le N) \in \reals^N: 
x_{(1)} \le x_{(2)} \le \cdots \le x_{(N)}  \right\}
\en
with constant drift vector $(\delta_j, 1 \le j \le N)$ and normal reflection at each of the
$N-1$ boundary hyperplanes $\{x_{(i)} = x_{(i+1)}\}$ for $1 \le i \le N-1$.
\end{lemma}

\proof
Sznitman \cite[page 594]{sznitman86}, \cite[Lemma 3.7]{sznitman91} gave these results in the
case of zero drifts.  As he observed, they follow from Tanaka's formula \cite[page 223]{ry99} and the definition 
of Brownian motion in a polyhedron with normal reflection. See, for example, Varadhan and Williams \cite{varwilliams}. Sznitman further proves that the Brownian motions $\{\beta_{j}\}$ are independent, which is essentially due to Knight's theorem \cite[page 183]{ry99}. 
The factors of $\sqrt{2}$ are most easily easily checked in the case $N=2$. That they are the same for all $N$ follows by a localization argument. The results of the lemma for general drifts $\delta_j$ are deduced from the results with no drift by application of the next two lemmas.
\endpf

\begin{lemma}\label{lemextra}
Let $X_1,X_2,\ldots$ be a sequence of independent Brownian motions indexed by a finite set $I$. Consider the ordered process $X_{(1)}, X_{(2)},\ldots$. Then, the following holds almost surely:
\eq\label{doublecollision}
\sum_{i=1}^{\abs{I}} \int_0^t 1\left\{\;X_i(s)=X_{(j)}(s)\;\right\}ds = t, \quad \forall\; j< \abs{I}+1,\; \forall\;t\in [0,\infty).
\en 
Moreover, the points of increments of the finite variation processes $L_{(j,j+1)}$ are almost surely disjoint.

When $I$ is countable, equation \eqref{doublecollision} still holds true as long as the assumptions on the initial values $X_1(0) < X_2(0) < \ldots$ are sufficient to guarantee the existence of the ordered processes for all finite times.
\end{lemma}

\proof For any two indices $k < l $, the lebesgue measure of the set $\{\;t\in [0,\infty):\; X_k=X_l\;\}$ is zero. This follows because the zero set of Brownian motion is of lebesgue measure zero almost surely. 

Now, the event $\sum_{i=1}^{\abs{I}} \int_0^t 1\left\{\;X_i(s)=X_{(j)}(s)\;\right\}ds > t$ for some $j$ and some $t$ implies that there exists some pair $(k,l)$ such that the lebesgue measure of the set $\{ 0\le s\le t:\; X_k(s)=X_{(j)}(s)=X_l(s)\;\}$ is positive. This is of measure zero according to the previous paragraph and by countable additivity.

The other possibility of $\sum_{i=1}^{\abs{I}} \int_0^t 1\left\{\;X_i(s)=X_{(j)}(s)\;\right\}ds < t$ is trivially ruled out by our assumption that the ordered processes are always achieved. 

For the second assertion, note that, according to the general theory of semi-martingale local times \cite{ry99}, the process $L_{(j,j+1)}$
increases only on the random closed set of times $t$ when $X_{(j)}(t) = X_{(j+1)}(t)$. These random sets are almost surely disjoint as $j$ varies. This is because, with probability $1$, there is {\em no triple collision}, meaning a time $t>0$ at which $X_i(t) = X_j(t) = X_k(t)$. It follows from the fact that the bivariate process $(X_i-X_j, X_j - X_k)$ is a linear transformation of a standard planar Brownian motion which does not hit points. 
\endpf

Call $\fsmpsp$ the {\em canonical setup} if $\fsmpsp$ is the usual space of continuous paths in $\reals^N$, with the usual right-continuous filtration, and $X_i$ is the $i$th coordinate process.

\begin{lemma}\label{lm2}
Let $\delta = (\delta_i, 1 \le i \le N) \in \reals^N$ and let $\mu$ be a probability distribution
on $\reals^N$.

\noindent{\em (i)} 
In the canonical setup there is a unique
probability measure $\PR^{\delta,\mu}$ under which the coordinate processes $(X_i, 1\le i\le N)$ solve the system
of SDE's \rem{sde1} with initial distribution $\mu$.
In particular, for $\delta =0$, the law $\PR^{0,\mu}$ is the Wiener measure governing 
standard Brownian motion in $\reals^N$ with initial distribution $\mu$.

\noindent{\em (ii)} 
In the canonical setup, for each $t >0$, the law $\PR^{\delta,\mu}$ is absolutely continuous
with respect to $\PR^{0,\mu}$ on $\FF_t$, with density
\eq
\lb{girs} 
\exp \left( \sum_{j = 1}^N \delta_j  \beta_{j}(t) - \hf \sum_{j=1}^N \delta_j ^2 t \right) 
\en
where $\beta_{j}$, which is the same as in \rem{tanak}, can be defined by the expression
\eq\label{whatisbetaj}
\beta_j(t) = \sum_{i=1}^N\int_0^t 1\left\{\; X_i(s)=X_{(j)}(s)\right\}dX_i(s),\qquad 1\le j\le N. 
\en
Under $\PR^{\delta,\mu}$ the $\beta_{j}$'s are independent Brownian motions on $\reals$ with drift coeffients $\{\delta_j\}$ and unit diffusion coefficients.

\noindent{\em (iii)} 
If $(X_i, 1 \le i \le N)$ is a realization of the $N$ particle system with drifts $\{\delta_i\}$
and initial distribution $\mu$ on an arbitrary probability space $(\Omega, \FF,\PR)$, then the 
$\PR$ joint distribution of the processes $X_i$ is identical to the $\PR^{\delta,\mu}$ distribution of 
the coordinate processes on the canonical space, as specified in {\em (ii)}. 

\end{lemma}
\proof
This is an instance of a well known general construction of the solution of an SDE with drift terms from one with no drift terms \cite[Chapter IX, Theorem (1.10)]{ry99}.

Under $\PR^{0,\mu}$, the fact that $\{\beta_j\}$ is a collection of independent Brownian motions also follows from expression \eqref{whatisbetaj} and equation \eqref{doublecollision}.
Under $\PR^{\delta,\mu}$, the process $\beta_j$ is a stochastic integral with respect to Brownian motions with drifts. By expanding the drift term, it is obvious that $\beta_j$ is a Brownian motion itself with drift $\delta_j$.
\endpf

Note that the SDE defining the $N$ particle system with drifts is a typical example of an SDE for which there is uniqueness in law, but not pathwise uniqueness.

For a solution $X$ of SDE \eqref{sde1}, let
$$
\bar{X}(t) := {1 \over N} \sum_{i = 1}^N X_i(t) = {1 \over N} \sum_{j = 1}^N X_{(j)}(t)
$$
which is the {\em center of mass} of the particle system, where we
regard each particle as having mass $1/N$.
We call the $\reals^N$-valued process
$$
(X_i - \bar{X}, 1 \le i \le N)
$$
the {\em centered system.} Note that the $N-1$ spacings defined by
taking differences of order statistics of the original system
are identical to the $N-1$ spacings defined by differences of 
the order statistics of the centered system.
So the $N$ order statistics of the centered system, which are constrained
to have average $0$, are an invertible linear transformation
of the $N-1$ spacings of original system.

\begin{lemma}\label{lm3}
For the $N$ particle system, with arbitrary drifts and initial distribution,

\noindent{\em (i)}
the center of mass process is a Brownian
motion with drift 
\eq
\lb{dave}
\bar{\delta}_N:=  N^{-1} \sum_{j = 1 }^N \delta_j
\en
and diffusion coefficient $1/N$; explicitly
\eq
\lb{Brep}
\bar{X}(t) = \bar{X}(0) + \bar{\delta} _N t + { B(t) \over \sqrt{N} }
\en
where 
$B(t):= N^{- \hf} \sum_{j = 1}^ N B_j(t)$ is a standard BM on $\reals$;
consequently
$$
\frac{\bar{X}(t)}{t} \te \bar{\delta}_N \mbox{ almost surely as}\; t \te \infty .
$$

\noindent{\em (ii)} The shifted center of mass process $\bar{X}(t)-\bar{X}(0)$ is independent of the centered system $\left(\;X_i(\cdot)- \bar{X}(\cdot)\;\right)$, $1\le i\le n$.

\end{lemma}
\proof Part (i) follows plainly from the SDE \re{sde1} by summing over $i$. 

For part (ii), note that if $\sigma_1$ and $\sigma_2$ are two independent sub-$\sigma$-algebras of a probability space $(A,\mcal{A},P)$, and we change $P$ to another probability $Q$ via defining $dQ/dP=fg$, where $f\in\sigma_1$, $g\in \sigma_2$, then $\sigma_1$ and $\sigma_2$ remain independent under $Q$. 

Now as in Lemma \ref{lm2}, when $(X_i)$ denotes $N$ independent Brownian motions with initial distribution $\mu$, the shifted average $\bar{X}(t)-\bar{X}(0)$ and the centered process $Y=X-\bar{X}\mathbf{1}$ are independent. This follows from the fact that conditionally on $X_0$, the processes are Gaussian with zero covariance, and that $\bar{X}(t)-\bar{X}(0)$ is independent of $X_0$. Now to get to $\PR^{\delta,\mu}$, the Radon-Nikodym derivative is given by \eqref{girs}. Note that, from expression \eqref{whatisbetaj} we get
\[ 
\beta_j(t)= \sum_{i=1}^N \int_0^t 1\left\{ Y_i(s) = Y_{(j)}(s) \right\}dY_i(s) + \bar{X}(t)-\bar{X}(0).
\]
Thus, from \eqref{girs}, it is clear that $d\PR^{\delta,\mu}/d\PR^{0,\mu}$ can be written as $fg$, where $f \in \sigma(Y)$ and $g \in \sigma(\bar{X}(\cdot) - \bar{X}(0))$. Now, by first localizing at finite time, and using the argument in the preceding paragraph, we establish the claim in part (ii).
\endpf

\begin{theorem}\label{theoremN} 
For $1 \le k \le N$ let
\eq
\lb{alphak}
\alpha_k:= \sum_{i = 1}^k ( \delta_i -  \bar{\delta}_N )  
\en
where $\bar{\de}_N$ is the average drift as in \rem{dave}.
For each fixed initial distribution of the $N$ particle system with
drifts $\{\delta_i, i=1,2,\ldots,N\}$, the collection of laws of 
$X_{(N)}(t) - X_{(1)}(t)$ for $t \ge 0$ is tight if
and only if 
\eq
\lb{conal}
\alpha _k >0  \mbox{ for all } 1 \le k \le N-1,
\en
in which case the following results all hold:

\noindent{\em (i)}
The distribution of the spacings system $(X_{(j+1)} - X_{(j)},\;1\le j\le N-1)$ at time $t$ converges in total
variation norm as $t \te \infty$ to a unique stationary distribution for
the spacings system, which is that of independent exponential
variables $Y_j$ with rates $2\alpha_j$, $1 \le j \le N-1$.
Moreover, the spacings system is reversible at equilibrium.

\noindent{\em (ii)}
The distribution of the centered system at time $t$ converges in total
variation norm as $t \te \infty$ to a unique stationary distribution for
the centered system, which is the distribution of
$$
( S_{\pi(i)-1} - \bar{S}, 1 \le i \le N )
$$
where $S_0 := 0$ and $S_i := Y_1 + \cdots + Y_i$ for $1\le i \le N-1$,
where $\pi$ a uniform random permutation of $\{1, \ldots, N\}$
which is independent of the $Y_i$, and
$$
\bar{S}:= {1 \over N} \sum_{i = 1}^N S_{\pi(i)-1} = {1 \over N} \sum_{i = 1}^{N-1} S_i = {1 \over N} \sum_{i = 1}^{N-1} (N-i) Y_i  .
$$
Moreover, the centered system is reversible at equilibrium.

\noindent{\em (iii)}
As $t \te \infty$
$$
X_i(t)/t \te \bar{\delta}_N \mbox{ almost surely for each } 1 \le i \le N ,
$$
and the same is true for $X_{(i)}(t)/t$ instead of $X_{i}(t)/t$.
\end{theorem}

\paragraph{Remark}
Regard the system as  split into {\em left-hand particles}
of rank $1$ to $k$ and {\em right-hand particles} of rank $k+1$ to $N$. 
If these two parts of the system are started at some strictly positive 
distance from each other, they evolve independently like copies of the 
$k$-particle system and the $(N-k)$-particle system respectively, until
the first time there is a collision between a left-hand particle
and a right-hand particle. It follows by summing up over the corresponding drifts that the centers of mass of the two subsystems left to themselves would have almost sure asymptotic speeds $\bar{\delta}_k$ and $\hat{\delta}_k$
respectively, where
$$
\bar{\delta}_k:= {1 \over k } \sum_{i = 1}^k \delta_i
\mbox{ and }
\hat{\delta}_k:= (N-k)^{-1} \sum_{i = k+1}^N \delta_i.
$$
Since
$$
\alpha_k = k \bar{\delta}_k - k \bar{\delta}_N  = { k (N-k) \over N } ( \bar{\delta}_k - \hat{\delta}_k )
$$
we see that $\alpha_k >0$ iff $\bar{\delta}_k > \hat{\delta}_k$,
which ensures that the righthand system cannot avoid an eventual
collision with the lefthand system.
According to \re{alphak},
for arbitrary prescribed $\bar{\delta}_N \in \reals$,
and $\alpha_k >0$, the unique drift vector determining an ergodic $N$-particle system whose average
drift is $\bar{\delta}_N$ and whose asymptotic spacings are independent exponential variables with rates $2 \alpha_k$, is given by
\eq
\lb{invs}
\delta_i = \bar{\delta}_N + \alpha_i - \alpha_{i-1} ~~~~~(1 \le i \le N)
\en
where $\alpha_0:= \alpha_N := 0$. 
Given an arbitrary cumulative probability distribution function $F$ on 
the line, 
and arbitrary $\bar{\delta} \in \reals$ and
$\eps >0$,
it is clear that by taking $N$ suitably large, we can choose
$(\alpha_k, 1 \le k \le N-1)$ and hence
$(\delta_j ,1 \le j \le N)$ so that for all sufficiently large $t$
$$
\PR \left( \sup_{x} \left| {1 \over N } \sum_{i = 1}^N 1( X_i (t) - \bar{X}(t) \le x )  - F(x)  \right| > \eps \right) < \eps 
$$
and $\bar{X}(t)/t \te \bar{\delta}$ almost surely as $t \te \infty$. Thus no 
matter what its initial distribution, the $N$ particle system looks 
asymptotically like a cloud of particles with mass distribution close
to $F$, which is drifting along the line at speed $\bar{\delta}$.

We would also like to mention here that a special case of the above theorem has been considered in a recent article by Jourdain and Malrieu \cite{joumal}. They consider SDE \eqref{sde1} with an increasing sequence of $\delta_i$'s and establishe joint convergence of the spacing system to independent exponentials as $t$ goes to infinity. 

\proof[Proof of Theorem \ref{theoremN}]
According to Lemma \ref{lm1}, the ordered $N$ particle system is a Brownian motion
$(X_{(k)}(t), 1 \le k \le N)_{t \ge 0}$ 
in the domain \re{domn} with identity covariance matrix, constant drift vector 
$(\delta_j, 1 \le j \le N)$ and normal reflections at each of the $N-1$ boundary hyperplanes 
$\{x_{(i)} = x_{(i+1)}\}$ for $1 \le i \le N-1$. Note that the vector $(1,1,\ldots,1)$ of all ones is in the intersection of these boundary hyperplanes. 

Let $\beta$ be an independent one-dimensional Brownian motion ($\beta_0=1$) with a negative drift $-\theta$ ($\theta >0$) which is reflected at the origin. The process defined by
\[
Y'_i(t) := X_i(t) - \bar{X}(t) + \frac{\beta}{\sqrt{N}},\quad i=1,2,\ldots N,
\]
has the same spacings as the original process $X$. Moreover, by the independence of the centered system and the center of mass process established in Lemma \ref{lm3}, and the fact that $\bar{Y'}=\beta/\sqrt{N}$ is the projection of $Y'$ on the subspace generated by the vector $(1,1,\ldots,1)$, it follows that the vector $(Y'_1,Y'_2,\ldots,Y'_N)$ is a Brownian motion with identity covariance matrix which is normally reflected in the wedge
\[
\left\{ y \in \rr^N: \sum_{i=1}^N y_i \ge 0, \; y_1 \le y_2\le\ldots\le y_N \right\}.
\]
The drift vector in this wedge can be written as a linear combination of spacings by summation by parts:
\eq
\lb{sumparts}
\sum_{i = 1}^N \left(\delta_i - \bar{\delta}_N -\frac{\theta}{\sqrt{N}}\right) y_i = 
- \sum_{k = 1}^{N-1} 
\left( \sum_{i = 1}^k (\delta_i - \bar{\delta}_N) \right)
(y_{k+1} - y_k) - \frac{\theta}{\sqrt{N}}\sum_{i=1}^N y_i .
\en
Noting that $Y'_{(k+1)}-Y'_{(k)}\equiv X_{(k+1)} - X_{(k)}$ for all $k$, the condition \re{conal} for stability, and part (i) of the Lemma, are now read from the general result about equilibrium distributions of reflecting Brownian motions stated in the following lemma, and standard theory of Harris recurrent Markov processes ( see, e.g., \cite{durrett}).

Part (ii) is established by showing that the centered system
$(X_i(t) - \bar{X}(t), 1 \le i \le N), t \ge 0$ is a Harris positive-recurrent diffusion
with the indicated invariant measure. The recurrence property has been proved in detail in \cite{atlasmodel}. The invariance is evident because a uniform randomization of labels relative to the
centered order statistics is clearly invariant for the centered motion. This convergence in distribution, combined with part (i) of Lemma \ref{lm2} gives convergence in probability of $X_i(t)/t $ to $\bar{\delta}_N$ for each $1 \le i \le N$. Almost sure convergence can now be justified by appeal to an ergodic theorem for the Harris recurrent centered diffusion.  
\endpf

The proof of Theorem \ref{theoremN} is completed by the following lemma, which we deduce from the general theory of stationary
distributions for reflecting Brownian motions in polyhedra due to
R. Williams \cite{williams87r}:
\begin{lemma}\label{lemmahw}
Let $R:= (R_t, t \ge 0)$ be a Brownian motion in the domain 
$$\{x \in \reals^K : b_i(x) \ge 0 \mbox{ for } i = 1,...K \}$$
for some collection of $K$ linearly independent linear functionals $b_i$,
with $R$ having identity covariance matrix, normal reflection at the boundary, and constant
drift $\delta$  with
\eq
\lb{db}
\sum_{i=1}^K \delta _i x_i = - \sum _{i =1}^K a_i b_i(x)  ~~~~~~(x_1, \ldots, x_K) = x \in \reals^K .
\en
This process $R$ has a stationary probability distribution iff $a_i >0 $ for all $i =1, \ldots, K$,
in which case in the stationary state the $b_i(R_t)$ are independent
exponential variables with rates $2 a_i$, and the process in its
stationary state is reversible.
\end{lemma}
\proof
This is read from the particular case of \cite[Theorem 1.2]{williams87r}
when the matrix $Q$ is identically $0$. 
According to that theorem, $R$ is in duality with itself relative to the
measure $\rho$ on the domain whose density function with respect to 
Lebesgue measure at $x$ is $\exp( 2 \delta \cdot x)$.
Using \re{db},
the linear change of variables to $y_i = b_i(x)$ shows that 
the $\rho$ distribution of the $b_i(x)$ has joint density at 
$(y_1, \ldots, y_K)$ equal to $c \prod_{i = 1}^K \exp( - 2 a_i y_i)$
for some $c >0$.
It follows that 
\eq
\lb{finmas}
\mbox{ $\rho$ has finite total mass iff $a_i >0$ for all $i$ }
\en
in which case, when $\rho$ is normalized to be a probability,
the 
$\rho$ distribution of the $b_i(x)$ is that of independent exponential
variables with rates $2 a_i$.
The ``if'' part of the conclusion is now evident. For the ``only if''
part we argue that if a stationary probability distribution $\rho'$
existed, it would obviously have a strictly positive density on the domain.
Then $R$ sampled at time $0,1,2,\ldots$ would be an irreducible Harris 
recurrent Markov process with respect $\rho'$, hence $\rho = c \rho'$ for 
some $c >0$ by the uniqueness of the invariant measure of a Harris recurrent 
Markov chain, and then $a_i >0$ for all $i$ by \re{finmas}.
\endpf
\comment{I guess that if some $a_i \le 0$ then $R$ is null-recurrent
if at most two of the $a_i$ equal $0$ and the rest are strictly
positive, and transient otherwise. Can we prove that?}

From Theorem \ref{theoremN} we immediately deduce:

\begin{crl}{crl3}
For each $\delta >0$ the $N$-particle Atlas system with drift vector 
$(\delta,0,\ldots,0)$
is ergodic with average speed $\delta/N$.
The stationary distribution of
$$
( X_{(j+1)} - X_{(j)}, 1 \le j \le N-1)
$$
is that of independent exponentials $(\zeta_j,\; 1\le j\le N-1)$, where the rate of $\zeta_j$ is $2\delta(1-j/N)$.
\end{crl}

\section{The infinite Atlas model}\label{sec.lag}

The infinite Atlas model can be described loosely as a countable collection of linear Brownian motions such that at every time point the minimum Brownian motion is given a positive drift of $\delta >0$, and the rest are left untouched. This is an example of \eqref{sde1} where $I=\mathbb{N}$, $\delta_1=\delta$ and all other $\delta_i=0$. Throughout this section we take $\delta=1$, since for our purpose here, the general case follows from the case when $\delta=1$ by scaling. 

It can be constructed rigorously in the weak solution framework in the following way. Start with the canonical setup (as in Lemma \ref{lm2}) of Brownian path-space:
\begin{equation}\label{pathspace0}
\left( \;\; C[0,\infty],\;\; \{\mcal{F}_t\}_{0\le t <\infty},\;\; \wien^x\;\;\right ),
\end{equation}
where $\{\mcal{F}_t\}$ is the right continuous filtration generated by the coordinates which satisfy the usual conditions, and $\wien^x$ is the law of the Brownian motion starting from $x\in \rr$. We now look at the countable product of a sequence of spaces like \eqref{pathspace0}:
\begin{equation}\label{pathspace}
\Omega = C[0,\infty]^{\mathbb{N}},\;\; \mcal{G}_t = \bigotimes_{1}^{\infty}\mcal{F}_t(i), \;\; P^x=\bigotimes_{1}^{\infty}\wien^{x_i},
\end{equation}
where the natural coordinate mapping is a countable collection of independent Brownian motions under $P^x$ starting from the sequence $x=(x_1,x_2,\ldots)$. 

Let $x$ be a sequence such that $x_{(1)} > -\infty$ and let $X=(X_1,X_2,\ldots)$ denote the sequence of infinite independent Brownian motions starting from $x$. We have the following lemma whose proof will follow later.

\begin{lemma}\label{infchange} Assume that the initial sequence $x$ is arranged in increasing order $x_1 \le x_2 \le \ldots$ and satisfies 
\begin{equation}\label{stability}
\liminf_{n \rightarrow \infty} \frac{(x_n - x_1)^2}{n} > 0.
\end{equation} 
Then, $P^x$-almost surely, $X_{(1)}(t) > -\infty$ for all $t \ge 0$, and the process
\begin{equation}\label{infmart}
N_t := \sum_{i=1}^{\infty} \int_0^t 1_{\{X_i(s) = X_{(1)}(s)\}}d X_i(s).
\end{equation}
is a $\{\mcal{G}_t\}$-martingale whose quadratic variation $\qvar{N}_t\equiv t$. The stochastic exponential of $N$, given by
\begin{equation}\label{whatisdt}
D_t = \exp\left(N_t - t/2\right).
\end{equation}
is hence again a non-negative $\{\mcal{G}_t\}$-martingale.
\end{lemma}

We change the measure $P^x$ by using the martingale $D$, i.e., define $Q$ by
\[
Q^x\vert_{\mcal{G}_t} \; := \; D_t\cdot P^x\vert_{\mcal{G}_t},\qquad t \ge 0.
\]
By Girsanov's theorem the probability measure $Q^x$ exists and is well-defined, and the coordinate process under $Q^x$ is a solution of the infinite Atlas model \eqref{sde1}. Hence $Q^x$ is the {\it{unique}} law of the Atlas model starting at $x$.
\bigskip

Our aim is the following: suppose the initial points $X_1(0) < X_2(0) < X_3(0) < \ldots$ are spread according to the Poisson process with rate two on the positive half-line, and we run the infinite Atlas model starting from these points. We shall prove that the product law of independent Exponential($2$) is invariant under the dynamics of the vector process $\Delta$ of spacings given by
\begin{equation}\label{whatisy}
\Delta_i(t) = X_{(i+1)}(t) - X_{(i)}(t),\qquad i=1,2,\ldots.
\end{equation}

The proof is achieved through a series of lemmas, the main argument being comparing the infinite Atlas model with the finite Atlas model and suitably passing to the limit. Throughout the proof the probability space is given by \eqref{pathspace}.

We start with the following lemma whose proof follows directly from Lemma \ref{lemextra}. We will find it convenient to use the following notation for the operator which sorts a given finite vector. If $x\in \rr^n,\; 1\le n < \infty$,  define
\begin{equation}\label{sortop}
\sort(x)= (x_{(1)},x_{(2)},\ldots,x_{(n)}).
\end{equation}

\begin{lemma}\label{coupling}
Every $\omega \in \Omega$ comprises of a sequence of processes $\omega(t) = ( \omega_1(t), \omega_2(t), \ldots )$. 
For any $N \in \mathbb{N}$, let us denote the ordered values of the processes with the first $N$ indices by
\[
Z^N(\omega) = (Z_1^N,Z_2^N,\ldots, Z_N^N)(\omega) = \sort((\omega_1, \omega_2, \ldots, \omega_N)).
\]
Then
\begin{equation}\label{whatisdnt}
D_N(t) = \exp\left(\sum_{i=1}^N \int_0^t 1_{\{X_i(s)=Z_1^N(s)\}}d X_i(s)  - t/2 \right)
\end{equation}
is a $\mcal{G}_t$-martingale. We denote by $Q^x_N$ the probability measure obtained by the changing the measure $P^x$ by the martingale $D_N$.
\end{lemma}

Since, under every $P^x$, each of the Brownian motions are independent, it follows, by applying Girsanov's Theorem, that, under $Q^x_N$, the first $N$ coordinates $(\omega_1, \ldots, \omega_N)$ evolve according to the finite Atlas model, while the rest of the coordinates are independent Brownian motions with the corresponding initial starting points. 
\medskip

Let $\mu$ denote the probability measure whereby the points $0=X_1(0) < X_2(0) < \ldots$ are such that the spacings $X_{i+1}(0) - X_i(0)$ are iid Exponential($2$).

\begin{lemma}\label{everythingexists}
For $\mu$-almost every $x$, the measure $Q^x$ exists and we can define
\[
\qmu = \int Q^x d\mu(x).
\]
\end{lemma}

\begin{proof}
The proof follows from Lemma \ref{infchange} and the law of large numbers.
\end{proof}

We now have our main theorem in this section which proves that $\mu$ is an invariant measure for the spacings of the infinite Atlas model.

\begin{theorem}\label{derivatzero}
For any $K \in \mathbb{N}$, and any function $F:\rr^K \rightarrow \rr$, which is smooth and has compact support, and for any time $t$, we have
\begin{equation}\label{zeroderiv} 
\e^{\qmu}\left[ F(\Delta_{1}, \Delta_{2},\ldots, \Delta_{K} )(t)\right]=\e^\mu\left[ F(\Delta_{1}, \Delta_{2},\ldots, \Delta_{K})(0)\right].
\end{equation}
Here, as defined in \eqref{whatisy}, $\Delta_i$ is the $i$th spacing $X_{(i+1)} - X_{(i)}$. 
\end{theorem}

\begin{proof}[Proof of Theorem \ref{thpois}]
Since the previous result holds for a class of functions which determines the marginal distributions of a sequence valued process, Theorem \ref{thpois} follows readily for $\delta=1$. The theorem for the infinite Atlas model with a general $\delta >0$ follows by scaling. 
\end{proof}

\begin{proof}[Proof of Lemma \ref{infchange}] Fix an arbitrary time $T > 0$, and consider $x$ as in the assumption of the lemma. We have the following claim. 

\begin{claim}\label{wellseparated}
\[
P^x\left( \;\omega:\;\exists \; M(\omega)\; \text{s.t.}\; \forall\; n \ge M,\;\; \inf_{0\le s\le T} (X_n(s) - X_{(1)}(s)) > 0 \; \right) = 1.
\]
\end{claim}

$\bullet$ We prove the above claim by establishing the following
\begin{equation}\label{spacedstart}
P^x\left( \;\omega:\;\exists \; M(\omega)\; \text{s.t.}\; \forall\; n \ge M,\;\; \inf_{0\le s\le T} (X_n(s) - X_{1}(s)) > 0 \; \right) = 1.
\end{equation}
This follows by defining the events
\begin{equation}\label{whatisai}
A_i := \left\{ \;\inf_{0\le s\le T} (X_i(s) - X_{1}(s)) \le 0\; \right\},\quad i=1,2,\ldots 
\end{equation}
Since $X_i- X_1$ is a Brownian motion starting from $(x_i - x_1)$, which for sufficiently large $i$ is strictly positive (by \eqref{stability}), by Bernstein's inequality \cite[page 153]{ry99}, one can easily estimate $P^x(A_i) \le \exp(-(x_i - x_1)^2/2T)$ for all large enough $i$. Again by assumption \eqref{stability}, we get
\[
\begin{split}
\limsup_{n \rightarrow \infty} \left( P^x(A_n) \right)^{1/n} &\le \limsup_n \exp(-(x_n - x_1)^2/2Tn)\\
& =  \exp\left(-\frac{1}{2T}\liminf_n\frac{(x_n - x_1)^2}{n}\right) < 1.\\
\end{split}
\]
Thus, by {Cauchy's root test}, it follows that the series sum $\sum_{i=1}^{\infty} P^x(A_i) < \infty$.
One can now apply the Borel-Cantelli lemma to obtain that $P^x(\limsup_i A_i)=0$, which proves \eqref{spacedstart}, and hence the required claim.
\medskip

We now prove that the process $N$ exists and is a martingale in time interval $[0,T]$. Since $T$ is arbitrary, this proves Lemma \ref{infchange}. Define the finite approximations
\[
B_k(t) = \sum_{i=1}^k \int_0^t 1_{\{X_i(s)=X_{(1)}(s)\}}dX_i(s), \quad 0 \le t \le T,\quad k=1,2,\ldots.
\]
Each $B_k$ is a stochastic integral with bounded, progressively measurable integrands. It is clear that they are martingales with quadratic variations
\[
\begin{split}
\qvar{B_k}_t &= \sum_{i=1}^k \int_0^t 1_{\{X_i(s)=X_{(1)}(s)\}}ds\\
&= \mathbb{L}\left\{\; 0 \le s\le t: \; \min_{1\le i \le k}X_i(s) = X_{(1)}(s)\;\right\},
\end{split}
\] 
where $\mathbb{L}$ refers to the Lebesgue measure on the line. We shall show that the sequence of martingales $\{B_k\}$ is a cauchy sequence in the $\mathbb{H}^2$-norm  , and hence has a limit which is denoted by $N$, as in \eqref{infmart}. 

To see this, observe that for any $n, k \in \mathbb{N}$, we get
\[
\qvar{B_{n+k} - B_n}_T = \mathbb{L}\left\{\; 0 \le s\le T: \; \min_{n+1\le i \le n+k}X_i(s) = X_{(1)}(s)\;\right\}.
\]
By Claim \ref{wellseparated}, we see that $P^x(\;\lim_{n,k \rightarrow \infty}\qvar{B_{n+k} - B_n}_T = 0 \;)=1$. It is also clear that $\qvar{B_{n+k} - B_n}_T\le T$, for all $n,k \in \mathbb{N}$. It follows from the Dominated Convergence Theorem that
\[
\lim_{n,k \rightarrow \infty} E^{P^x}\left(\; \qvar{B_{n+k} - B_n}_T\;\right) =0. 
\]
This, by definition, shows that the sequence of martingales $\{B_n\}$ is a cauchy sequence in the $\mathbb{H}^2$ norm. Since the space of continuous martingales under that norm is complete, the sequence $\{B_n\}$ converges to a limiting martingale which we denote $N$. It also follows that $\qvar{N}_t = \lim_{n \rightarrow \infty} \qvar{B_n}_t = t.$ Thus, $N$ is actually a $\{\mcal{G}_t\}$-Brownian motion. This completes the proof of Lemma \ref{infchange}.
\end{proof}

\noindent{\bf Remark.} The condition in \eqref{stability} is clearly loose. In fact, if we consider $A_i$ as in \eqref{whatisai}, all we require is that, for every $T > 0$, one should have
\[
\sum_{i=1}^{\infty}P^x(A_i) \le \sum_{i=1}^{\infty}\exp\left(-(x_i-x_1)^2/2T\right) < \infty,
\]
which is a much weaker condition than required by \eqref{stability}.
\medskip

The proof of Theorem \ref{derivatzero} relies on probability estimates proved in the next three Lemmas \ref{KeyEstimate1}, \ref{KeyEstimate2}, and \ref{KeyEstimate3}. The second one is actually a generalization of the first. But, we choose to do the first separately since it is simpler and more transparent. But first we will do some basic computations.

\begin{lemma}\label{gammasq}
If $Y \sim Gamma(r,\lambda)$ for some $r \ge 1$, then, for any $t > 0$, we have 
\[
\begin{split}
\e\left(e^{-Y^2/2t}\right) &\le e^{t\lambda^2/2}\left( 2\lambda^2t\right)^{r/2}\frac{\Gamma(r/2)}{2\Gamma(r)}\\
&= \sqrt{\pi} \left(\frac{\lambda^2t}{2}\right)^{r/2}\frac{e^{t\lambda^2/2}}{\Gamma((r+1)/2)}.
\end{split}
\]
\end{lemma}

\begin{proof}
\eq\label{eqnran1}
\begin{split}
\e&\left(e^{-Y^2/2t}\right) = \int_0^{\infty}e^{-y^2/(2t)}\frac{\lambda^r}{\Gamma(r)}y^{r-1}e^{-\lambda y}dy\\
&= \frac{\lambda^r}{\Gamma(r)}e^{t\lambda^2/2}\int^{\infty}_0y^{r-1}e^{-(y+t\lambda)^2/(2t)}dy\\
&= \frac{\lambda^r}{\Gamma(r)}e^{t\lambda^2/2}\int^{\infty}_{\sqrt{t}\lambda}(\sqrt{t}z-t\lambda)^{r-1}e^{-z^2/2}\sqrt{t}dz, \quad z= \frac{y + t\lambda}{\sqrt{t}}.
\end{split}
\en
For $z \ge \sqrt{t}\lambda$, since $r\ge 1$, one has
\[
(\sqrt{t}z - t\lambda)^{r-1} \le (\sqrt{t}z)^{r-1}.
\]
Thus one can bound the \eqref{eqnran1} by  
\[
\begin{split}
\e\left(e^{-Y^2/2t}\right)&\le \frac{\lambda^r}{\Gamma(r)}e^{t\lambda^2/2}\int^{\infty}_{\sqrt{t}\lambda}(\sqrt{t}z)^{r-1}e^{-z^2/2}\sqrt{t}dz\\
&= \frac{(\lambda^2t)^{r/2}}{\Gamma(r)}e^{t\lambda^2/2}\int^{\infty}_{\sqrt{t}\lambda}z^{r-1}e^{-z^2/2}dz\\ 
&=\frac{(\lambda^2t)^{r/2}}{\Gamma(r)}e^{t\lambda^2/2}\int^{\infty}_{t\lambda^2/2}(2w)^{(r-1)/2}e^{-w}(2w)^{-1/2}dw, \quad w=z^2/2\\
&=\frac{(\lambda^2t)^{r/2}}{\Gamma(r)}e^{t\lambda^2/2}2^{r/2-1}\int_{t\lambda^2/2}^{\infty}w^{r/2-1}e^{-w}dw\\
&\le \frac{(2\lambda^2t)^{r/2}}{2\Gamma(r)}e^{t\lambda^2/2}\Gamma(r/2).  
\end{split}
\]
The final identity in the lemma is due to the \emph{duplication formula}:
\[
\frac{\Gamma(r/2)}{\Gamma(r)}= \sqrt{\pi}\frac{2^{1-r}}{\Gamma((r+1)/2)}.
\]
\end{proof}

\begin{lemma}\label{KeyEstimate1}
For any $t > 0$, and all positive integers $N$ satisfying $N +1 \ge 16et$, under $\pmu$, we have the following estimate of the probability of the event that during time $[0,t]$, the globally lowest ranked process is in fact the lowest ranked process among the processes with the first $N$ indices.
\[
\pmu\left\{ X_{(1)}(s) = Z^N_1(s),\quad 0\le s\le t \right\} \ge 1- C_1 e^{2t}\left(\sqrt{\frac{4et}{N+1}} \right)^{N},
\]
where $C_1$ is a positive constant. To remind the reader, the processes $(Z^N_k, \; 1\le k \le N)$ have been defined in Lemma \ref{coupling}.
\end{lemma}

\begin{proof}
Note that the complement of the event has the following upper bound.
\begin{equation}\label{minapprox}
\begin{split}
1 &- \pmu\left\{ X_{(1)}(s) = Z^N_1(s),\quad 0\le s\le t \right\}\\
&\le \pmu\left(\;  \cup_{i\ge N+1}  \{\;X_i(s) \le Z^N_1(s), \;\text{for some}\; s\in [0,t]\;\}\; \right)\\
&\le \pmu\left(\;  \cup_{i\ge N+1}  \{\;X_i(s) \le X_1(s), \;\text{for some}\; s\in [0,t]\;\}\; \right)\\
&\le \sum_{i=N+1}^{\infty}\pmu\left(\; X_i(s) \le X_1(s), \;\text{for some}\; s\in [0,t] \;\right).
\end{split}
\end{equation}
The final bound above is the so-called union bound.

Now we will use exponential bounds for Brownian supremums to estimate $\pmu\left(\; X_i(s) \le X_1(s), \;\text{for some}\; s\in [0,t] \;\right)$ which is the same as $\pmu\left(\; \inf_{0\le s\le t}(X_i(s) - X_1(s)) \le 0 \;\right)$. Note that under $P\cdot\mu$, the process $X_i-X_1$ is a Brownian motion whose initial distribution is the law of $X_i(0)- X_1(0)$. This law is Gamma($i-1, 2$), since it is the sum of $(i-1)$ iid Exp($2$). Conditional on $X(0)=x$, such that $\{X_i(0)- X_1(0)=y\}$, it follows from Berstein's inequality that
\[
P^x\left(\; \inf_{0\le s\le t}(X_i(s) - X_1(s)) \le 0 \;\right) \le \exp(-y^2/2t).
\]
Thus, for $i \ge 2$, we define $Y=X_i(0)-X_1(0) \sim \text{Gamma}(i-1,2)$, and we use Lemma \ref{gammasq} to get
\begin{equation}\label{infbound}
\begin{split}
\pmu\left(\; \inf_{0\le s\le t}(X_i(s) - X_1(s)) \le 0 \;\right) &\le \e(\exp(-Y^2/2t))\\
&\le \sqrt{\pi}(2t)^{(i-1)/2}\frac{e^{2t}}{\Gamma(i/2)}.
\end{split}
\end{equation}

Plugging the estimate in \eqref{minapprox}, we get
\begin{equation}\label{minapprox2}
\begin{split}
1 - \pmu\left\{\; X_{(1)} = Z^N_1(s), \; 0\le s\le t \right\} &\le \sqrt{\pi}e^{2t}\sum_{i=N+1}^{\infty}\frac{(2t)^{(i-1)/2}}{\Gamma(i/2)}.
\end{split}
\end{equation}
By Stirling's approximation, $\exists\; C$, a positive constant, such that
\begin{equation}\label{stapp}
\Gamma(z)\ge C^{-1}e^{-z}z^{z-1/2}, \quad \forall\; z\in \rr^+.
\end{equation} 
Thus, from \eqref{minapprox2} we get
\[
\begin{split}
1 &- \pmu\left\{\; X_{(1)} = Z^N_1(s), \; 0\le s\le t \right\} \le C \sqrt{\pi}e^{2t}\sum_{i=N+1}^{\infty}\frac{(2t)^{(i-1)/2}e^{i/2}}{(i/2)^{(i-1)/2}}\\
&= C\sqrt{\pi}e^{2t}\sum_{i=N}^{\infty}\left(\frac{4et}{i+1}\right)^{i/2}e^{1/2}\le \sqrt{\pi}Ce^{2t+1/2} \sum_{i=N}^{\infty}\left(\sqrt{\frac{4et}{N+1}} \right)^{i}\\
&\le 2C\sqrt{\pi} e^{2t+1/2}\left(\sqrt{\frac{4et}{N+1}} \right)^{N},\quad \text{when}\; \sqrt{4et/(N+1)} \le 1/2.
\end{split}
\]
This proves the estimate.
\end{proof}

\begin{lemma}\label{KeyEstimate2}
For any any $t > 0$ and all positive integers $1\le k < N \in \mathbb{N}$ satisfying $N-k+2 \ge 16et$, under $\pmu$, we have the following estimate of the probability of the event that in the time interval $[0,t]$, the globally lowest $k$ ranked processes are in fact the lowest $k$ ranked processes among the ones with the first $N$ indices.
\begin{equation}\label{kest2}
\begin{split}
\pmu&\left\{ \left(X_{(1)}, \ldots, X_{(k)}\right)(s) = \left(             Z^N_1, \ldots, Z^N_k\right)(s), \quad 0\le s\le t \right\}\\ &\ge 1 -  C_2(k) e^{2t} (N-1)^{k-1} \left(\sqrt{\frac{4et}{N-k+2}}\right)^{N-k+1},
\end{split}
\end{equation}
for some positive constant $C_2$ depending on $k$. For $k=1$, we get back the bound in Lemma \ref{KeyEstimate1}.
\end{lemma}

\begin{proof}
As in the previous lemma, we bound the probability of complement of the event.
\begin{equation}\label{kapprox}
\begin{split}
1 &- \pmu\left\{ \left(X_{(1)}, \ldots, X_{(k)}\right)(s) = \left(             Z^N_1, \ldots, Z^N_k\right)(s), \quad 0\le s\le t \right\} \\
&\le \pmu\left(\;\cup_{i\ge N+1}\left\{\;X_i(s) \le Z^N_k(s),\quad\text{for some}\; s\in [0,t] \;   \right\}   \right)\\
&\le \sum_{i=N+1}^{\infty}\pmu\left\{\;X_i(s) \le Z^N_k(s),\quad\text{for some}\; s\in [0,t] \;   \right\}.
\end{split}
\end{equation}
Now, if $[N]$ denotes the set $\{1,2,\ldots,N\}$, let us note that
\[
Z_k^N = \max_{i_1 < i_2 < \ldots < i_{k-1}}\;\;\min_{l \in [N] \backslash \{i_1, i_2, \ldots, i_{k-1}\}} X_l.
\]
Thus, to get an upper bound on $\pmu\left\{\;X_i(s) \le Z^N_k(s),\quad\text{for some}\; s\in [0,t] \;   \right\}$, one can apply the union bound one more time to obtain 
\begin{multline}\label{kapprox2}
\pmu\left\{\;X_i(s) \le Z^N_k(s),\quad\text{for some}\; s\in [0,t] \;   \right\}\le \\
\sum_{i_1 <  \ldots < i_{k-1}}\pmu\left( X_i(s) \le \min_{l \in [N] \backslash \{i_1, \ldots, i_{k-1}\}} X_l(s),\quad\text{for some}\; s\in [0,t] \;\right) \le \\
\sum_{i_1 < \ldots < i_{k-1}}\pmu\left( X_i(s) \le X_{i^*}, \text{for some}\; s\in [0,t], i^*=\min \{ [N] \backslash \{i_1,\ldots, i_{k-1}\}\} \;\right)
\end{multline}
Now, we can count the frequency with which $i^*$ takes its possible values as the choice $i_1 < i_2 < \ldots < i_{k-1}$ varies in $\{ 1, 2, \ldots, N\}$. Let $g(i)$ be the number of ways to pick $i_1 < i_2 < \ldots < i_{k-1}$ such that $i^*=i$. Then, it is straightforward to see
\[
\begin{split}
g(1) &= \#\{\;\{i\}:\; i_1 > 1 \;\}=\combi{N-1}{k-1},\\
g(2) &= \#\{\;\{i\}:\; i_1=1,\; i_2 > 2\;\}= \combi{N-2}{k-2}\\
&\vdots\\
g(l) &= \#\{\;\{i\}:\; i_1=1,i_2=2,\ldots,i_{l-1}=l-1,\; i_l > l\;\}= \combi{N-l}{k-l},\;\; l \le k,\\
g(l)&=0,\quad \forall\; l > k.
\end{split}
\]
Thus, by \eqref{kapprox2}, we get
\begin{multline}\label{kapprox3}
\sum_{i_1 < i_2 < \ldots < i_{k-1}}\pmu\left( X_i(s) \le \min_{l \in [N] \backslash \{i_1,\ldots, i_{k-1}\}} X_l(s),\quad\text{for some}\; s\in [0,t] \;\right)\\
\le \sum_{l=1}^{k}\combi{N-l}{k-l}\pmu\left( X_i(s) \le X_{l}(s),\;\text{for some}\; s\in [0,t]\;\right).
\end{multline}
Now, again as we did when deriving the last key estimate, $X_i-X_l$ for $i > l$ is a Brownian motion under $\pmu$ with the initial distribution Gamma($i-l,2$). Thus, by Lemma \ref{gammasq}, we can bound
\[
\pmu\left( X_i(s) \le X_{l}(s),\;\text{for some}\; s\in [0,t]\;\right) \le \sqrt{\pi}\frac{\left(2t\right)^{(i-l)/2}e^{2t}}{\Gamma\left((i-l+1)/2\right)}.
\]
Combining the above inequality with estimates \eqref{kapprox}, \eqref{kapprox2}, and \eqref{kapprox3}, we get
\begin{equation}\label{kapprox4}
\begin{split}
1 &- \pmu\left\{ \left(X_{(1)}, \ldots, X_{(k)}\right)(s) = \left(             Z^N_1, \ldots, Z^N_k\right)(s), \quad 0\le s\le t \right\} \\
&\le \sum_{i=N+1}^{\infty} \sum_{l=1}^k \combi{N-l}{k-l}\sqrt{\pi}e^{2t}\frac{\left(2t\right)^{(i-l)/2}}{\Gamma\left((i-l+1)/2\right)}.
\end{split}
\end{equation}
We will again give a loose but good enough upper bound for the infinite sum. But, first we need to note that, for any $n \ge k$, we have
\[
{\combi{n}{k}}\Big/{\combi{n-1}{k-1}}= \frac{n}{k} \ge 1.
\]
Thus, it follows that
\[
\combi{N-l}{k-l} \le \combi{N-1}{k-1},\quad \forall\; 1\le l\le k.
\]
We can use this to simplify \eqref{kapprox4}:
\[
\begin{split}
1 &- \pmu\left\{ \left(X_{(1)}, \ldots, X_{(k)}\right)(s) = \left(             Z^N_1, \ldots, Z^N_k\right)(s), \quad 0\le s\le t \right\} \\
&\le \sum_{i=N+1}^{\infty} \sum_{l=1}^k \combi{N-1}{k-1}\sqrt{\pi}e^{2t}\frac{\left(2t\right)^{(i-l)/2}}{\Gamma\left((i-l+1)/2\right)}\\
\end{split}
\]
We again use Stirling's approximation \eqref{stapp} to get
\begin{equation}\label{kapprox5}
\begin{split}
\sum_{i=N+1}^{\infty}& \sum_{l=1}^k \combi{N-1}{k-1}\sqrt{\pi}e^{2t}\frac{\left(2t\right)^{(i-l)/2}}{\Gamma\left((i-l+1)/2\right)}\\
&\le C\combi{N-1}{k-1}\sqrt{\pi}e^{2t}\sum_{i=N+1}^{\infty}\sum_{l=1}^k \frac{(2t)^{(i-l)/2}e^{(i-l+1)/2}}{((i-l+1)/2)^{(i-l)/2}}\\
&= C\combi{N-1}{k-1}\sqrt{\pi}e^{2t+1/2}\sum_{i=N+1}^{\infty}\sum_{l=1}^k \left(\sqrt{\frac{4et}{i-l+1}}\right)^{i-l}.
\end{split}
\end{equation}
Since $i-l \ge N+1-k$ for $l \le k < N+1 \le i$, it is clear that $4et/(i-l+1) \le 4et/(N+2-k)$. And hence
\[
\sum_{l=1}^k \left(\sqrt{\frac{4et}{i-l+1}} \right)^{i-l} \le \sum_{l=1}^k \left(\sqrt{\frac{4et}{N-k+2}} \right)^{i-l}.
\]
By our assumption
\begin{equation}\label{ass2}
\sqrt{\frac{4et}{N-k+2}}\le \frac{1}{2}
\end{equation}
Thus
\[
\left(\sqrt{\frac{4et}{N-k+2}} \right)^{i-l} \le \left(\sqrt{\frac{4et}{N-k+2}} \right)^{i-k},\quad l \le k < i,
\]
and consequently
\[
\sum_{l=1}^k \left(\sqrt{\frac{4et}{i-l+1}} \right)^{i-l} \le k\left(\sqrt{\frac{4et}{N-k+2}} \right)^{i-k}.
\]
Plugging this bound in \eqref{kapprox5} we get
\begin{equation}\label{kapprox51}
\begin{split}
&\sum_{i=N+1}^{\infty} \sum_{l=1}^k \combi{N-1}{k-1}\sqrt{\pi}e^{2t}\frac{\left(2t\right)^{(i-l)/2}}{\Gamma\left((i-l+1)/2\right)}\\
&\le C\combi{N-1}{k-1}\sqrt{\pi}e^{2t+1/2}\sum_{i=N+1}^{\infty}  k\left(\sqrt{\frac{4et}{N-k+2}} \right)^{i-k}\\
&\le Ck\combi{N-1}{k-1}\sqrt{\pi}e^{2t+1/2}\left(\sqrt{\frac{4et}{N-k+2}}\right)^{N+1-k}\sum_{j=0}^{\infty} \left(\sqrt{\frac{4et}{N-k+2}}\right)^{j}\\
&\le C_2(k) e^{2t} (N-1)^{k-1} \left(\sqrt{\frac{4et}{N-k+2}}\right)^{N-k+1},\; \text{by}\; \eqref{ass2}, 
\end{split}
\end{equation}
for some positive constant $C_2$ depending on $k$. This proves the lemma.
\end{proof}

\begin{lemma}\label{KeyEstimate3} For any $N \in \mathbb{N}$, define $\mu_N$ to be the law under which 
\begin{equation}\label{whatismun}
\begin{split}
X_1(0)=0,\quad X_{i+1}(0) - X_i(0) &\sim \text{Exp}\left(2( 1 - i/N) \right), \quad i=1,2,\ldots, N-1,\\
\text{and}\quad X_{i+1}(0) - X_i(0) &\sim \text{Exp}\left(2 \right), \quad i=N,N+1,\ldots,
\end{split}
\end{equation}
and all these spacings are independent.

For any $t > 0$, and any three integers $k < J < N$ satisfying
\begin{equation}\label{ass3}
J-k+2 \ge 16et,
\end{equation} 
we have the following bound on the probability that under $\pmu_N$, during the time interval $[0,t]$, the lowest $k$ ranked processes among the processes with the first $N$ indices are in fact the lowest $k$ among those with the first $J$ indices. 
\begin{equation}\label{kest3}
\begin{split}
&\pmu_N\left\{\; \left(Z^J_{1}, \ldots, Z^J_{k}\right)(s) = \left(             Z^N_1, \ldots, Z^N_k\right)(s), \quad 0\le s\le t\; \right\}\\ 
&\ge 1-  C_2(k)e^{2t}(J-1)^{k-1}\left(\sqrt{\frac{4et}{J-k+2}}\right)^{J-k+1}\left[ 1  - \left(\sqrt{\frac{4et}{J-k+2}}\right)^{N-J}\right].
\end{split}
\end{equation}
\end{lemma}

\begin{proof}
We follow the same line of argument as in the last two estimates. Thus
\[
\begin{split}
1 &- \pmu_N\left\{\;\left(Z^J_{1}, \ldots, Z^J_{k}\right)(s) = \left(             Z^N_1, \ldots, Z^N_k\right)(s), \quad 0\le s\le t \;\right\}\\
&\le \sum_{i=J+1}^N \pmu_N\left(\; X_i(s)\le Z^J_k(s),\quad\text{for some}\;s\in [0,t] \;\right). 
\end{split}
\]
Following similar counting arguments as in \eqref{kapprox2} and \eqref{kapprox3}, we can bound
\begin{equation}\label{kapproxn2}
\begin{split}
1 &- \pmu_N\left\{\;\left(Z^J_{1}, \ldots, Z^J_{k}\right)(s) = \left(             Z^N_1, \ldots, Z^N_k\right)(s), \quad 0\le s\le t \;\right\}\\
&\le \sum_{i=J+1}^N \sum_{l=1}^k \combi{J-l}{k-l}\pmu_N\left(\; X_i(s)\le X_l(s),\quad\text{for some}\;s\in [0,t] \;\right). 
\end{split}
\end{equation}

Now, under $\mu_N$, the gap $X_i(0)-X_l(0)=\sum_{j=l}^{i-1}Y_j$, where the $Y_j$'s are independent and $Y_j$ is distributed as Exp($2(1-j/N)$). Thus, by the exponential bound used before, we get
\begin{equation}\label{kapproxn1}
\begin{split}
\pmu_N&\left(\; X_i(s)\le X_l(s),\quad\text{for some}\;s\in [0,t] \;\right)\\
&\le \e\left[ \exp\left(-\frac{1}{2t}\left(\sum_{j=l}^{i-1}Y_j\right)^2 \right)\right].
\end{split}
\end{equation}
Now, since each $Y_j$ is Exponential($2(1-j/N)$), the random variables 
\[
Y_j^*= (1-j/N)Y_j
\]
are iid Exp($2$), and each $Y_j^* \le Y_j$ since $j \le N$. Thus, $\sum_{j=l}^{i-1}Y_j^* \le \sum_{j=l}^{i-1}Y_j$, and hence
\[
\e\left[ \exp\left(-\frac{1}{2t}\left(\sum_{j=l}^{i-1}Y_j\right)^2 \right)\right] \le \e\left[ \exp\left(-\frac{1}{2t}\left(\sum_{j=l}^{i-1}Y^*_j\right)^2 \right)\right].
\]
But, each $Y_j^*$ is Exp($2$), and hence
\[
\e\left[ \exp\left(-\frac{1}{2t}\left(\sum_{j=l}^{i-1}Y^*_j\right)^2 \right)\right] = \e\left[e^{-Y^2/2t}\right],
\]
where $Y$ is a Gamma($i-l,2$) random variable. Thus, plugging in the bound from Lemma \ref{gammasq} into \eqref{kapproxn1}, we derive
\[
\begin{split}
\pmu_N\left(\; X_i(s)\le X_l(s),\quad\text{for some}\;s\in [0,t] \;\right) &\le  \sqrt{\pi}\frac{e^{2t}(2t)^{(i-l)/2}}{\Gamma((i-l+1)/2)}.
\end{split}
\]
Now, we follow approximations similar to \eqref{kapprox5} and \eqref{kapprox51} for the right-hand-side on \eqref{kapproxn2} to obtain
\[
\begin{split}
&1 - \pmu_N\left\{\;\left(Z^J_{1}, \ldots, Z^J_{k}\right)(s) = \left(             Z^N_1, \ldots, Z^N_k\right)(s), \quad 0\le s\le t \;\right\}\\
&\le \sum_{i=J+1}^N \sum_{l=1}^k \combi{J-1}{k-1}\sqrt{\pi}\frac{e^{2t}(2t)^{(i-l)/2}}{\Gamma((i-l+1)/2)}\\
&=C\combi{J-1}{k-1}\sqrt{\pi}e^{2t+1/2}\sum_{i=J+1}^N \sum_{l=1}^k \left(\sqrt{\frac{4et}{(i-l+1)}}\right)^{i-l}\\
&\le C\combi{J-1}{k-1}\sqrt{\pi}e^{2t+1/2}\sum_{i=J+1}^N k \left(\sqrt{\frac{4et}{(J-k+2)}}\right)^{i-k}, \; \text{by}\; \eqref{ass3},\\
&\le C_2(k)e^{2t}(J-1)^{k-1}\left(\sqrt{\frac{4et}{J-k+2}}\right)^{J-k+1}\quad\sum_{j=0}^{N-J-1}\left( \sqrt{\frac{4et}{J-k+2}} \right)^j.
\end{split}
\]
Now, if we call $r=\sqrt{4et/(J-k+2)}$, then $r\le 1/2$ by asusmption \eqref{ass3}. The finite geometric sum can be easily bounded as
\[
\sum_{j=0}^{N-J-1}r^j = \frac{1-r^{N-J}}{1-r}\le 2(1-r^{N-J}).
\]
Thus, suitably altering the constant $C_2$, we get our desired bound
\[
\begin{split}
1 &- \pmu_N\left\{\;\left(Z^J_{1}, \ldots, Z^J_{k}\right)(s) = \left(             Z^N_1, \ldots, Z^N_k\right)(s), \quad 0\le s\le t \;\right\}\\
&\le C_2(k)e^{2t}(J-1)^{k-1}\left(\sqrt{\frac{4et}{J-k+2}}\right)^{J-k+1}\left[ 1  - \left(\sqrt{\frac{4et}{J-k+2}}\right)^{N-J}\right].
\end{split}
\]
This proves the lemma.
\end{proof}

\bigskip

\begin{proof}[Proof of Theorem \ref{derivatzero}] For every $K\le N\in\mathbb{N}$, let
\[
Y_i^N=Z_{i+1}^N - Z_i^N, \quad 1\le i\le N.
\]
From Corollary \ref{crl3}, we know that the law of the first $(N-1)$ spacings under $\mu_N$ (defined in Lemma \ref{KeyEstimate3}) is exactly the stationary distribution of spacings for the finite Atlas model with $N$ particles. Since, under $Q_N$ (see Lemma \ref{coupling}), the dynamics of the processes $\{X_1(t), \ldots,X_N(t)\}$ is that of a finite Atlas model independent of the rest of the Brownian motions, it follows by stationarity that for any $t > 0$, we have 
\begin{equation}\label{finstat}
\e^{\qmun{N}}\left[ F(Y_1^N, \ldots, Y_K^N)(t)\right]= \e^{\mu_N}\left[F(Y_1^N, \ldots, Y_K^N)(0)\right].
\end{equation}
We will show that for fixed $t$ and $K$, as $N$ tends to infinity, the two sides of the above equation converges to the corresponding sides of \eqref{zeroderiv}. This will prove the theorem.

However, to do this we will need an intermediary stage where for an integer $J < N$, the dynamics of the process is according to $Q_J$ while the initial distribution of the spacings is either $\mu$ or $\mu_N$. Define, for $K < J < N$, the following quantities:
\begin{eqnarray*}
a&=&\e^{\qmun{N}}\left[ F(Y_1^N, \ldots, Y_K^N)(t)\right]-\e^{\qjmun{J}{N}}\left[ F(Y_1^J, \ldots, Y_K^J)(t)\right],\\
b&=&\e^{\qjmun{J}{N}}\left[ F(Y_1^J,  \ldots, Y_K^J)(t)\right]-\e^{\qjmu{J}}\left[ F(Y_1^J, \ldots, Y_K^J)(t)\right],\\
c&=&\e^{\qjmu{J}}\left[ F(Y_1^J, \ldots, Y_K^J)(t)\right]- \e^{\qmu}\left[ F(\Delta_1, \ldots, \Delta_K)(t)\right],\;\text{and}\\
d&=&\e^{\mu}\left[ F(\Delta_1, \ldots, \Delta_K)(0)\right]- \e^{\mu_N}\left[F(Y_1^N, \ldots, Y_K^N)(0)\right].
\end{eqnarray*}
It is clear that $a,b,c$, and $d$, all depend on $t,K,J,$ and $N$, although we choose to suppress this dependence in the notation. Also, it follows plainly from their definitions combined with equality \eqref{finstat} that 
\[
\left\lvert\; \e^{\qmu}\left[ F(\Delta_1, \ldots, \Delta_K)(t)\right] - \e^{\mu}\left[ F(\Delta_1, \ldots, \Delta_K)(0)\right]\;\right\rvert\;\le\; \abs{a} + \abs{b} + \abs{c} + \abs{d}.
\]
We will now show that $a,b,c,d$ all go to zero if we select a sequence of $J$ and $N$ such that $J$, $N$ and $N/J^2$ go to infinity. This will prove the theorem. 
\medskip

\noindent{\bf Step 1.}[Estimate of $a$.] For $x_1\le x_2 \le \ldots \le x_{K+1}$, define 
\[
G(x_1,x_2,\ldots,x_{K+1}):=F(x_2-x_1,x_3-x_2,\ldots,x_{K+1}-x_K).
\] 
Then, clearly, $G$ is also a continuous bounded function. Assume that
$\sup_{x}\abs{G(x)} \le \alpha$. Now, by the changing the measures from $Q_N$ and $Q_J$ to $P$, we get
\begin{equation}\label{smalla1}
\begin{split}
a &= \e^{\qmun{N}}\left[G(Z_1^N,\ldots,Z_{K+1}^N)(t)\right] - \e^{\qjmun{J}{N}}\left[G(Z_1^J,\ldots,Z_{K+1}^J)(t)\right]\\
&= \e^{\pmu_{N}}\left[D_N(t)G(Z_1^N,\ldots,Z_{K+1}^N)(t)- D_J(t)G(Z_1^J,\ldots,Z_{K+1}^J)(t)\right],
\end{split}
\end{equation}
where $D_N$ and $D_J$ are the Radon-Nikod\'ym derivative processes defined in \eqref{whatisdnt}. Define the event 
\[
\Gamma := \left\{\; \left(Z_1^N,\ldots,Z_{K+1}^N\right)(s)= \left(Z_1^J,\ldots,Z_{K+1}^J\right)(s),\; \forall\; s\in [0,t]\;\right\}.
\]
If $\omega\in \Gamma$, by definition, $\pmu_N$ almost surely $D_N(t,\omega)=D_J(t,\omega)$. Thus \eqref{smalla1} can be written as
\[
\begin{split}
a&= \e^{\pmu_{N}}\left[\;\left(D_N(t)G(Z_1^N,\ldots,Z_{K+1}^N)- D_J(t)G(Z_1^J,\ldots,Z_{K+1}^J)\right)1_{\Gamma^c}\;\right]
\end{split}
\]
Thus, $\abs{a}$ is bounded above by
\begin{equation}\label{smalla2}
\begin{split}
&\e^{\pmu_{N}}\left\lvert\;D_N(t)G(\cdot)1_{\Gamma^c}\right\rvert + \e^{\pmu_{N}}\left\lvert\;D_J(t)G(\cdot)1_{\Gamma^c}\;\right\rvert\\
&\le \alpha \left( \e^{\pmu_{N}}\left\lvert\;D_N(t)1_{\Gamma^c}\right\rvert + \e^{\pmu_{N}}\left\lvert\;D_J(t)1_{\Gamma^c}\;\right\rvert   \right)\\
&\le \alpha \left(\norm{D_N(t)}_N + \norm{D_J(t)}_N\right)\; \sqrt{\pmu_N\left( \Gamma^c\right)}.
\end{split}
\end{equation}
The final inequality is due to the Cauchy-Schwartz inequality, where the norm $\norm{\cdot}_N$ refers to the $\ltwo$ norm under the measure $\pmu_N$. 

Now, under $\pmu_N$, both the Radon-Nikodym derivatives $D_N(t)$ and $D_J(t)$ are equal in law to $\exp(B_t - t/2)$, where $B$ is a standard Brownian motion. Thus, it is straightforward to see
\[
\norm{D_N(t)}_N = \norm{D_J(t)}_N = \exp(t/2).
\]
Also, by Lemma \ref{KeyEstimate3} (put $k=K+1$ in the statement), for large enough $J$ and $N$ we have $\pmu_N\left(\Gamma^c\right)$ is less than
\[
C_2(K+1)e^{2t}(J-1)^{K}\left(\sqrt{\frac{4et}{J-K+1}}\right)^{J-K}\left[ 1  - \left(\frac{4et}{J-K+1}\right)^{N-J}\right].
\]
If we plug everything back to \eqref{smalla2}, we see that $\abs{a}$ goes to zero as $J$, $N$, and $N/J^2$ go to infinity, while keeping $t$ and $K$ fixed.
\medskip

\noindent{\bf Step 2.}[Estimate of $b$.] Under $Q^x_J$, the vector $(Y_1^J,\ldots,Y_K^J)$ depends only on the first $J$ processes $X_1,X_2,\ldots,X_J$, and independent of $X_i,\;i> J$. Thus
\[
\e^{Q^x_J}\left[ F(Y_1^J,\ldots,Y_K^J)\right] = H(x_1,x_2,\ldots,x_J),
\]
where $H$ is bounded function since $F$ is bounded. Thus, we have the following equality
\[
\begin{split}
\abs{b}&= \abs{\;\e^{\mu_N}(H(X_1(0),\ldots,X_J(0))) - \e^{\mu}(H(X_1(0),\ldots,X_J(0)))\;}
\end{split}
\]
If $\nu_N$ and $\nu$ denote the law of the vector $(X_1(0),\ldots,X_J(0))$ under $\mu_N$ and $\mu$ respectively, then, since we have assume $G$ (or equivalently $F$) to be bounded in absolute value by $\alpha$, we have $\abs{b}\le \alpha \norm{\nu_N-\nu}_{TV}$. Here, $\norm{\cdot}_{TV}$ refers to the total variation norm. We will now show that $\norm{\nu_N-\nu}_{TV}$ goes to zero as $N$, $J$, and $N/J^2$ go to infinity.

Under $\mu_N$, $X_1(0)=0$ and each initial spacing $X_{i+1}(0)-X_{i}(0)$, $1\le i \le N-1$, follows independent Exp($2(1-i/N)$). While, under $\mu$, $X_1(0)=0$ and the spacings follow iid Exp($2$). Now the law of the vector $(X_i(0),\;1\le i\le J)$ is determined by the first $J$ spacings $(X_{i+1}(0)- X_i(0), \; 1\le i\le J)$ which gives us the following inequality 
\[
\begin{split}
\norm{\nu_N - \nu}_{TV}&\le \int_{\rr^J}\left\lvert \prod_{i=1}^J 2(1-i/N)e^{-2(1-i/N)x_i} - 2^Je^{-2\sum_{i=1}^J x_i} \right\rvert dx\\
&=\int_{\rr^J}\left\lvert \prod_{i=1}^J (1-i/N)e^{2ix_i/N} - 1 \right\rvert2^Je^{-2\sum_{i=1}^J x_i} dx.
\end{split}
\]
By an application of Cauchy-Schwartz inequality, we get
\begin{equation}\label{boundtv}
\begin{split}
\norm{\nu_N-\nu}_{TV}^2 &\le  \int_{\rr^J} \left\lvert \prod_{i=1}^J(1-i/N)e^{2ix_i/N} - 1 \right\rvert^2 2^Je^{-2\sum_{1}^Jx_i}dx,\\
&=2^J\prod_{i=1}^J(1-i/N)^2  \int_{\rr^J}\exp\left\{{-2}\sum_{i=1}^J (1-2i/N)x_i\right\}dx\\
&- 2^{J+1}\prod_{i=1}^J(1-i/N) \int_{\rr^J}\exp\left\{{-2}\sum_{i=1}^J (1-i/N)x_i\right\}dx + 1
\end{split}
\end{equation}
By the standard identity $\int e^{-\lambda x}dx = \lambda^{-1}$ for $\lambda > 0$, we get
\eq\label{bbound}
\begin{split}
\norm{\nu_N-\nu}_{TV}^2&\le 2^J\prod_{i=1}^J(1-i/N)^2 \prod_{i=1}^J (2(1-2i/N))^{-1}\\
& - 2^{J+1}\prod_{i=1}^J(1-i/N)\prod_{i=1}^J (2(1-i/N))^{-1} + 1\\
&= \frac{(1-1/N)^2(1-2/N)^2\ldots(1-J/N)^2}{(1-2/N)(1-4/N)\ldots(1-2J/N)} - 1.
\end{split}
\en

The following inequality is straightforward to prove
\[
e^{-2x} \le 1- x \le e^{-x}, \quad \text{for all}\quad 0 \le x \le \frac{1}{2}\log 2.
\]
By our assumption $J$, $N$, and $N/J^2$ are going to infinity, we can assume for all sufficiently large values of $J$ and $N$ that $2J/N \le \log 2/2$. By the previous inequality, we get
\[
e^{-2i/N}\le (1-i/N) \le e^{-i/N}, \quad 1\le i \le 2J/N,
\]
and consequently
\begin{multline}\label{bbound2}
\frac{(1-1/N)^2(1-2/N)^2\ldots(1-J/N)^2}{(1-2/N)(1-4/N)\ldots(1-2J/N)}\le \frac{\exp\left( -2\sum_{i=1}^J i/N \right)}{\exp\left( -2\sum_{i=1}^J 2i/N \right)}\\
\text{and}\quad\frac{\exp\left( -4\sum_{i=1}^J i/N \right)}{\exp\left( -\sum_{i=1}^J 2i/N \right)}\le \frac{(1-1/N)^2(1-2/N)^2\ldots(1-J/N)^2}{(1-2/N)(1-4/N)\ldots(1-2J/N)}.
\end{multline}
Thus, we can give upper and lower bounds
\[
e^{-{J(J+1)}/{N}} \le\frac{(1-1/N)^2(1-2/N)^2\ldots(1-J/N)^2}{(1-2/N)(1-4/N)\ldots(1-2J/N)}\le e^{{J(J+1)}/{N}}.
\]
Combining this with \eqref{bbound} we see that as $N$, $J$, and $N/J^2$ go to infinity, we clearly get that $\abs{b}$ converges to zero.
\medskip

\noindent{\textbf{Step 3.}}[Estimate of $c$.]
This is similar to the methods we use to estimate $a$. As in {Step 1}, we first employ a change of measure to obtain 
\[
\begin{split}
c &= \e^{\pmu}[D_{J}(t)G(Z_1^J,\ldots,Z_{K+1}^J)(t)] - \e^{\pmu}[D(t)G(X_{(1)}, \ldots, X_{(K+1)})(t)].
\end{split}
\]
Consider the event
\begin{equation}\label{whatisgamma}
\begin{split}
\Gamma(t,J,K)&=\left\{\left(X_{(1)}, \ldots, X_{(K+1)}\right)(s) = \left( Z^J_1, \ldots, Z^J_{K+1}\right)(s), \; 0\le s\le t \right\}
\end{split}
\end{equation}
On $\Gamma(t, J,K)$, the random processes $G(Z_1^J,\ldots,Z_{K+1}^J)(s)$ are identical to $G(X_{(1)}, \ldots, X_{(K+1)})(s)$ in time $s\in[0,t]$. Also, the processes $Z^J_1$ and $X_{(1)}$ are the same in the time interval $[0,t]$. Thus the process $D(t)$ and $D_J(t)$ are also the same. Since, by our assumption, $\abs{G}$ is bounded by $\alpha$, the following upper bound on $\abs{c}$ holds
\begin{equation}\label{cissmall}
\begin{split}
\abs{c} &\le \alpha\e^{\pmu}\left[ \left(\abs{D(t)} +\abs{D_{J}(t)}\right)1_{\Gamma^c(t,J,K)}\right].\\
&\le \alpha\left(\norm{D_{t}} + \norm{D_{J}(t)}  \right)\sqrt{\pmu(\Gamma^c(t,J,K))}, 
\end{split}
\end{equation}
where, we denote by $\norm{\cdot}$ the $\ltwo$ norm under the measure $\pmu$.

Now, by Lemma \ref{KeyEstimate2} (for $k=K+1$ and $N=J$ in the statement), for large enough $J$ we get
\begin{equation}\label{boundaepsilon}
\begin{split}
\pmu(\;\Gamma^c(t,J,K)\;) &\le C_2(K+1) e^{2t} (J-1)^{K} \left(\sqrt{\frac{4et}{J-K+1}}\right)^{J-K}. 
\end{split}
\end{equation}
Now, since $t$ and $K$ are fixed, as $J$ tends to infinity, $\pmu(\;\Gamma^c(t,J,K)\;)$ goes to zero. Finally, as in the estimate of $a$ in Step 1, note that, by \eqref{whatisdt} and $\eqref{whatisdnt}$, we can assert that  
\[
D(t)= \exp(X_{t}-t/2),\quad\text{and}\quad D_N({t})= \exp(Y_{t}-t/2),
\]
where $X$ and $Y$ are Brownian motions. Thus $\norm{D(t)}=e^{t/2}=\norm{D_N(t)}$.
If we plus back everything in \eqref{cissmall}, we get $\abs{c}$ goes to zero as $J$ tends to infinity.
\medskip

\noindent{\bf Step 4.}[Estimate of $d$.] At time zero  the indices are arranged in the increasing order, and hence $(\Delta_1,\Delta_2,\ldots, \Delta_K)$ is obviously equal to $(Y_1^N,Y_2^N,\ldots, Y_K^N)$. Since $\abs{F}\le \alpha$, it follows that $\abs{d}$ is bounded by $\alpha$ times the total variation distance between the law of $(\Delta_1,\Delta_2,\ldots, \Delta_K)$ under $\mu$ and $\mu_N$. A computation exactly similar to \eqref{boundtv} gives us
\[
\left(\frac{\abs{d}}{\alpha}\right)^2 \le \frac{(1-1/N)^2(1-2/N)^2\ldots(1-K/N)^2}{(1-2/N)(1-4/N)\ldots(1-2K/N)} - 1.
\]
But, since $K$ is fixed, and $N$ grows to infinity, the right-hand-side above goes to zero by a logic similar to \eqref{bbound2}. This proves the estimate.

Thus, combining Steps $1,2,3$ and $4$, we have proved the theorem.
\end{proof}

\section*{Acknowledgements}
We thank Prof. M. E. Vares for drawing our attention to articles \cite{demasiferrari} and \cite{rostvares}. We are also grateful to the anonymous referee for an extremely helpful review which significantly improved the quality of the paper, including pointing out references \cite{dgl85} and \cite{dgl87}.

\def\cprime{$'$} \def\polhk#1{\setbox0=\hbox{#1}{\ooalign{\hidewidth
  \lower1.5ex\hbox{`}\hidewidth\crcr\unhbox0}}} \def\cprime{$'$}
  \def\cprime{$'$} \def\cprime{$'$}
  \def\polhk#1{\setbox0=\hbox{#1}{\ooalign{\hidewidth
  \lower1.5ex\hbox{`}\hidewidth\crcr\unhbox0}}} \def\cprime{$'$}
  \def\cprime{$'$} \def\polhk#1{\setbox0=\hbox{#1}{\ooalign{\hidewidth
  \lower1.5ex\hbox{`}\hidewidth\crcr\unhbox0}}} \def\cprime{$'$}
  \def\cprime{$'$} \def\cydot{\leavevmode\raise.4ex\hbox{.}} \def\cprime{$'$}
  \def\cprime{$'$} \def\cprime{$'$}



\end{document}